\documentclass[11pt]{article}

\newtheorem{thm}{Theorem}[section]
\newtheorem{prop}[thm]{Proposition}
\newtheorem{lem}[thm]{Lemma}
\newtheorem{cor}[thm]{Corollary}

\newtheorem{example}[thm]{Example}
\newtheorem{remark}[thm]{Remark}

\newcommand{\pr}[1]{{I\!\!P}^{#1}}
\newcommand{\skipit}[1]{{}}
\newcommand{\prfend}{\hbox to7pt{\hfil}
\par\vskip-\baselineskip\hbox to\hsize
{\hfil\vbox {\hrule width6pt height6pt}}\vskip\baselineskip}

\newcommand{\myarrow}[2]{\hbox to #1pt{\hfil$\to$\hfil}{\hskip-#1pt{\raise
10pt\hbox to#1pt{\hfil$\scriptscriptstyle #2$\hfil}}}}

\begin{document}

\title{Resolutions
of ideals of fat points with support in a hyperplane}

\author{Giuliana Fatabbi\\
Dip. di Matematica e Informatica\\
Universit\`a di Perugia\\
via Vanvitelli 1\\
06123 Perugia, Italy\\
email: fatabbi@dipmat.unipg.it
\and
Brian Harbourne\\
Department of Mathematics\\
University of Nebraska\\
Lincoln, NE 68588-0130 USA\\
email: bharbour@math.unl.edu
\and
Anna Lorenzini\\
Dip. di Matematica e Informatica\\
Universit\`a di Perugia\\
via Vanvitelli 1\\
06123 Perugia, Italy\\
email: annalor@dipmat.unipg.it}

\maketitle

\thanks{Acknowledgments: We thank MURST,
whose national project {\it Algebra
Commutativa e Computazionale}, and
the University of
Perugia, whose project {\it Metodi
algebrici e analitici nello studio delle variet\`a}, supported
visits to Perugia by the second author, who also
thanks the NSA and NSF for supporting his research.
We also thank the referee for helpful suggestions.}

\begin{abstract}
Let $Z'$ be a fat point subscheme of $\pr{d}$, and let $x_0$ be
a linear form such that some power of $x_0$ vanishes on $Z'$
(i.e., the support of $Z'$ lies in the hyperplane $H$ defined by
$x_0=0$, regarded as $\pr{d-1}$). Let $Z(i)=H\cap Z'(i)$,
where $Z'(i)$ is the subscheme of $\pr{d}$ residual to $x_0^i$;
note that $Z(i)$ is a fat points subscheme of $\pr{d-1}=H$.
In this paper we give a graded free resolution of the ideal $I(Z')$ over
$R'=K[\pr{d}]$, in terms of the graded minimal free resolutions
of the ideals $I(Z(i))\subset R=K[\pr{d-1}]$. We also give
a criterion for when the resolution is minimal, and we show that
this criterion always holds if $\hbox{char}(K)=0$.
\end{abstract}

\section{Introduction}\label{intro}

Let $R=K[\pr {d-1}]=K[x_1,\ldots,x_d]$ and
$R'=K[\pr {d}]=K[x_0, x_1,\ldots,x_d]$ be the homogeneous
coordinate rings of projective space, over an algebraically closed
field $K$ of arbitrary characteristic. We regard $\pr {d-1}$
as the hyperplane $x_0=0$ in $\pr {d}$. We will denote
homogeneous components by subscripts; thus, for example, $R_1$
denotes the $K$-vector space of homogeneous linear forms in
$K[\pr {d-1}]$.

Given points $p_1,\ldots,p_r\in \pr {d-1}$ and
nonnegative integers $m_i$, we have the fat point subschemes
$Z=m_1p_1+\cdots+m_rp_r\subset\pr {d-1}$ (so $I(Z)\subset R$)
and $Z'=m_1p_1+\cdots+m_rp_r\subset \pr {d}$ (so
$I(Z')\subset R'$).
In particular, the ideal $I(Z)$ is $I_1^{m_1}\cap \cdots\cap I_r^{m_r}$,
where $I_j$ is the ideal generated by all forms in $R$
vanishing at $p_j$, and $I(Z')$ is $(I')_1^{m_1}\cap \cdots\cap (I')_r^{m_r}$,
where $(I')_j$ is the ideal generated by all forms in $R'$
vanishing at $p_j$.

We also have the obvious canonical inclusion
$R\subset R'$, so we can regard ideals in $R$ as $R$-submodules of
$R'$. Now define
$Z'_{m-i}=(m_1-i)_+p_1+\cdots+(m_r-i)_+p_r\subset \pr {d}$
and $Z_{m-i}=(m_1-i)_+p_1+\cdots+(m_r-i)_+p_r\subset
\pr {d-1}$ for each $0\le i\le
m=\hbox{max}\{m_1,\ldots,m_r\}$, where for any integer $n$ we
define $n_+=\hbox{max}\{n,0\}$. Note that $Z'=Z'_m$ and $Z=Z_m$,
and that $\emptyset=Z_0\subset Z_1\subset\cdots\subset Z_m=Z$.
Alternatively, let $Z'(i)$ denote the subscheme of $\pr{d}$ defined
by the ideal $I(Z'):(x_0^i)$ (and thus residual to $x_0^i$);
then $Z'(i)=Z'_{m-i}$ and $Z_{m-i}$ is the subscheme
$Z'(i)\cap \pr{d-1}$ of $\pr{d-1}$.

In this paper, we construct a graded free resolution of $I(Z')$ over $R'$,
given graded minimal free resolutions of each $I(Z_i)$ over $R$.
Under the condition that $I(Z_{i+1})\subset R_1I(Z_i)$
for each $i$, we show that the constructed resolution is minimal.
We also show that this condition always holds if $\hbox{char}(K)=0$.
In fact, we do not know any examples where the condition does not hold.

As a corollary we obtain a result about the Poincar\'e polynomial
of $I(Z')$. (Recall that the Poincar\'e polynomial encodes the
Betti numbers of a resolution. Given a subscheme $W\subset \pr n$
and its ideal $I(W)\subset A=K[\pr n]$,
we have its minimal free resolution $0\to F_t\to \cdots\to F_1\to
F_0\to I(W)\to 0$. The Poincar\'e
polynomial $P(W)$ is defined as follows. Each syzygy module $F_i$
is a free graded $A$-module of the form $F_j = \bigoplus_i
A(-i)^{a_{ij}}$, and we take
$$P(W) = \Sigma_{ij} a_{ij}T^iX^j\in {\bf Z}[X,T].$$
So, for example, the empty subscheme $W=\emptyset$, whose ideal is
thus $(1)$, has polynomial $P(W) = 1$, and if $W=p$ is a single
reduced point in $\pr 3$, then $P(W)=3T+3XT^2+X^2T^3$.)

\begin{thm}\label{Pointhm} Let $Z'=m_1p_1+\cdots+m_rp_r\subset \pr {d}$,
where the points $p_i$ lie in a hyperplane of $\pr {d}$.
If $Z_i$ is defined as above and $I(Z_{i+1})\subset R_1I(Z_i)$ for each $i$, then
$P(Z') = (1 + XT)(\Sigma_{0<i\le m} T^{m-i}P(Z_i)) + T^m$, where $m$ is the
maximum of the multiplicities $m_1,\ldots,m_r$.
\end{thm}

We give our resolution construction in Section \ref{constr}.
To most usefully apply our construction, we also need
examples in which the condition $I(Z_{i+1})\subset R_1I(Z_i)$
holds. More generally, given $p_i\in \pr d$ and fat points
$Z=m_1'p_1+\ldots+m_r'p_r$ and $Y=m_1p_1+\ldots+m_rp_r$ with
$m_i'>m_i$ for all $i$ whenever $m_i'>0$, it seems possible that
$I(Z)\subset R_1I(Y)$ always holds. We do not know of any counterexamples,
and we show in Section \ref{appls} various
situations, such as whenever $\hbox{char}(K)=0$, where this condition does hold.
We also discuss various examples in Section \ref{appls}. In particular, our
construction allows us
to determine the resolution of the ideal of two fat points in any
projective space, recovering some results (but not the explicit formulas)
of \cite {fatlor} and \cite {valla}. Since it is
easy to give the resolution of any number of fat points
in $\pr1$, our result also allows us, among other cases,
to determine the resolution for any number
of fat points in any projective space if the points lie
on a line.

One can speculate on possible generalizations.
In one direction, giving a resolution of an ideal $I\subset R'$
is equivalent to giving a resolution for $R'/I$, and in the cases
we study here the annihilator of $R'/I$ contains a power of a linear form.
Perhaps our approach can be applied
in other cases of modules over a polynomial ring annihilated
by a power of a linear form. But even for cyclic modules it is unclear
what kind of criterion to expect for resolutions to be minimal.
Since our immediate interest is minimal resolutions of ideals of fat
points, we have not pursued such questions here.
In another direction, one can ask to relax the condition of
considering points that lie only on a hyperplane.
However, this would seem to change the problem in a fundamental way,
since we implicitly use the fact that $R'$ modulo a linear form is
not only a quotient of $R'$ but also a subring of $R'$.

\section{The Construction}\label{constr}

We will need to refer to the following two elementary results.

\begin{lem}\label{fundlem}
Let $h: M\to N$ be a (not necessarily graded) homomorphism of
graded $R$-modules. Let $F$ and $G$ be free modules over $R$ with
surjective $R$-homomorphisms $\alpha:F\to M$ and $\beta:G\to N$.
Then there is an $R$-homomorphism $h_0:F\to G$ such that $\beta
h_0=h\alpha$, where $h_0$ is homogeneous of degree $t$ if $h$ is
homogeneous of degree $t.$ If furthermore $h(M)\subset R_1N$, then
in addition $h_0$ can be chosen such that $h_0(F)\subset R_1G$.
\end{lem}
{\em Proof.} The first part is clear, since $F$ is free, so assume
$h(M)\subset R_1N$ and consider $B:G^d\to R_1N\subset N$, where
$B:(c_1,\ldots,c_d)\mapsto x_1\beta(c_1)+\cdots+x_d\beta(c_d)$.
Clearly, $B$ is surjective, so we can lift $h\alpha:F\to M\to
R_1N$ to $h': F\to G^d$. We also have the canonical map
$\gamma:G^d\to G$ (in which $(c_1,\ldots,c_d)\mapsto
x_1c_1+\cdots+x_dc_d$), and $\beta\gamma=B$, and hence
$\beta\gamma h'=Bh'$. If we take $h_0=\gamma h'$, then $\beta
h_0=Bh'=h\alpha$ as desired, and $h_0(F)\subset \gamma(G^d)=R_1G$.
\prfend

By recursively applying the previous lemma, we
obtain:

\begin{cor}\label{fundcor}
Let $h: M\to N$ be a (not necessarily graded) homomorphism of
graded $R$-modules, with $F_\bullet$ and $G_\bullet$ free
resolutions over $R$ of $M$ and $N$ respectively. Then there are
$R$-homomorphisms $h_j:F_j\to G_j$, $j\ge0$, compatible with the
differential morphisms of the resolutions, where each $h_j$ is
homogeneous of degree $t$ if $h$ is homogeneous of degree $t$.
Moreover, if $h(M)\subset R_1N$, then the maps $h_j$ can be chosen
such that $h_j(F_j)\subset R_1G_j$ for every $j$.
\end{cor}

Now consider $\emptyset=Z_0\subset\cdots\subset Z_m= Z\subset \pr
{d-1}$ and $Z'\subset \pr {d}$ as in the introduction. Given
minimal graded free resolutions (over $R$) for each $I(Z_i)$, we
now construct a graded resolution (over $R'$) for $I(Z')$. We will
use the following notation: the graded free modules in the
resolution of $I(Z_i)$ will be denoted $F_{i,j}$ (so $F_{i,0}$ is
the free $R$-module on the generators of $I(Z_i)$, with the
suitable shifts; $F_{i,1}$ the free $R$-module on the first
syzygies of $I(Z_i)$, with the suitable shifts; etc.). The free
generators for $F_{i,j}$ will be denoted $s_{k,i,j}$, indexed by
$k$. The graded resolution differential $F_{i,j+1}\to F_{i,j}$
will be denoted $\phi_{i,j+1}$. We denote the  augmentation maps
$F_{i,0}\to I(Z_i)$ by $\phi_{i,0}$.

We will also need the maps $f_{i+1,j}:F_{i+1,j}\to F_{i,j}$
guaranteed by Corollary \ref{fundcor}, where, in the notation of
the corollary, $h:M\to N$ is the inclusion
$M=I(Z_{i+1})\subset I(Z_{i})=N$, and $f_{i+1,j}=h_j$. By the
corollary, each $f_{i+1,j}$ is homogeneous of degree 0 and
$f_{i+1,j}\phi_{i+1,j+1} = \phi_{i,j+1}f_{i+1,j+1}$, for $j\ge0,$
and $\phi_{i+1,0} = \phi_{i,0}f_{i+1,0}.$

We now construct a resolution of $I(Z')$ of the form
$$\cdots\myarrow{25}{\phi'_{j+2}}
F'_{j+1}\myarrow{25}{\phi'_{j+1}} F'_{j}\myarrow{25}{\phi'_{j}}
\cdots\myarrow{25}{\phi'_1} F'_0\myarrow{25}{\phi'_0}I(Z')\to 0.$$

To avoid repeatedly having to indicate certain shifts explicitly,
we denote $F_{i,j}(-(m-i))\otimes_R R'$ by $F_{i,j}'.$ Each map
$\phi_{i,j}:F_{i,j}\to F_{i,j-1}$ induces an obvious map
$F_{i,j}(-(m-i))\to F_{i,j-1}(-(m-i))$ which extends to give a map
$\phi'_{i,j}:F_{i,j}'\to F_{i,j-1}'$; i.e., after accounting for
the shift, $\phi'_{i,j}$ is just
$\phi_{i,j}\otimes_R\hbox{id}_{R'}$. Similarly, $f'_{i+1,j}:
F_{i+1,j}'(-1)\to F_{i,j}'$ denotes the map coming from
$f_{i+1,j}$. Now take
$$F'_0=\bigoplus_{i=0}^m F_{i,0}'\hbox to1.5in{\hfil and, for $j\ge1$,\hfil}
F'_j=F'_{0,j}\oplus\left(\bigoplus_{i=1}^m \left(F_{i,j}'\oplus
F_{i,j-1}'(-1)\right)\right)$$ 
(note $F'_{0,j}=0$ for $j\ge1$; we include it in $F'_j$ for consistency). Define the  augmentation map as
$$\phi'_0(s_{k,i,0}\otimes1)=\phi'_{i,0}(s_{k,i,0}\otimes
x_0^{m-i})$$
and differentials as
$$\phi'_1(s_{k,i,0}\otimes1)=s_{k,i,0}\otimes x_0 -
f'_{i,0}(s_{k,i,0}\otimes1),$$ and, for $j\ge1$,
$$\phi'_j(s_{k,i,j}\otimes1)=\phi'_{i,j}(s_{k,i,j}\otimes1),$$
and
$$\phi'_{j+1}(s_{k,i,j}\otimes1)=s_{k,i,j}\otimes x_0 -
f'_{i,j}(s_{k,i,j}\otimes1)-\phi'_{i,j}(s_{k,i,j}\otimes1).$$

Note the ambiguity of whether an element $s_{k,i,j-1}\otimes1$
lies in $F'_j$ or $F'_{j-1}$. We will resolve this ambiguity
either by an explicit declaration, such as $s_{k,i,j-1}\otimes1\in F'_j$,
or implicitly, as in $\phi'_j(s_{k,i,j-1}\otimes1)$ (keeping in mind that the
differential $\phi'_j$ is a mapping defined on $F'_j$).

The modules $F_j'$ and maps $\phi_j'$  can also be described  in terms of mapping cones.
In fact, we have the mapping
of complexes $f_{i+1}':F'_{i+1,\bullet}(-1)\to F'_{i,\bullet}$,
where for each $j$ the map $F'_{i+1,j}(-1)\to F'_{i,j}$ is given
by $f'_{i+1,j}$. Thus this mapping of complexes is ultimately
induced by the inclusion $Z_i\subset Z_{i+1}$. We also have the
mapping of complexes $\mu_{i}:F_{i,\bullet}'(-1)\to
F_{i,\bullet}'$, where $\mu_{i}(a)=x_0a$ is given by multiplying
by $x_0$. Then, for $j\ge0,$ the modules $F_j'$ and maps $\phi_{j+1}'$
can be seen as coming from amalgamating the mapping cones of
$\mu_{m}$, $f_{m}'$, $\cdots$, $\mu_{1}$, $f_{1}'$, as
shown in Figure 1. Any two consecutive rows
give a mapping cone, either for some $\mu$ or for some $f'$.
The direct sum of the modules in column  $j$ gives $F'_j$. The
differential $\phi'_j:F'_j\to F'_{j-1}$ is the direct sum of the
maps between columns. 

$$\begin{array}{cccccccc}
 & & & & & & & \\
\to & F_{m,3}' & {\phi_{m,3}'\atop\displaystyle\to} & F_{m,2}' & {\phi_{m,2}'\atop\displaystyle\to}
& F_{m,1}' & {\phi_{m,1}'\atop\displaystyle\to} & F_{m,0}'   \\
 & & & & & & & \\
     &  & \hbox to0in{\hss$\scriptstyle\mu_{m}$}\nearrow  & & \hbox to0in{\hss$\scriptstyle\mu_{m}$}\nearrow
     &  & \hbox to0in{\hss$\scriptstyle\mu_{m}$}\nearrow &  \\
 & & & & & & & \\
\to & F_{m,2}'(-1) & {-\phi_{m,2}'\atop\displaystyle\to}  & F_{m,1}'(-1)
& {-\phi_{m,1}'\atop\displaystyle\to}  & F_{m,0}'(-1) &  &  \\
 & & & & & & & \\
     &                &  \hbox to0in{\hss$\scriptstyle -f'_{m,2}$}\searrow &
     & \hbox to0in{\hss$\scriptstyle -f'_{m,1}$} \searrow &     &  \hbox to0in{\hss$\scriptstyle -f'_{m,0}$}\searrow &   \\
 & & & & & & & \\
\to & F_{m-1,3}' & {\phi_{m-1,3}'\atop\displaystyle\to} & F_{m-1,2}' & {\phi_{m-1,2}'\atop\displaystyle\to}
& F_{m-1,1}' & {\phi_{m-1,1}'\atop\displaystyle\to} & F_{m-1,0}'   \\
 & & & & & & & \\
     &  & \hbox to0in{\hss$\scriptstyle\mu_{m-1}$}\nearrow  &
     & \hbox to0in{\hss$\scriptstyle\mu_{m-1}$}\nearrow &  & \hbox to0in{\hss$\scriptstyle\mu_{m-1}$}\nearrow &  \\
 & & & & & & & \\
\to & F_{m-1,2}'(-1) & {-\phi_{m-1,2}'\atop\displaystyle\to}  & F_{m-1,1}'(-1)
& {-\phi_{m-1,1}'\atop\displaystyle\to}  & F_{m-1,0}'(-1) &  &  \\
 & & & & & & & \\
     &                &  \hbox to0in{\hss$\scriptstyle -f'_{m-1,2}$}\searrow &
     & \hbox to0in{\hss$\scriptstyle -f'_{m-1,1}$} \searrow &     &  \hbox to0in{\hss$\scriptstyle -f'_{m-1,0}$}\searrow &   \\
 & & & & & & & \\
     &                &  \cdots &   & \cdots &     & \cdots &   \\
 & & & & & & & \\
      &  & \hbox to0in{\hss$\scriptstyle\mu_{1}$}\nearrow  &
     & \hbox to0in{\hss$\scriptstyle\mu_{1}$}\nearrow &  & \hbox to0in{\hss$\scriptstyle\mu_{1}$}\nearrow &  \\
 & & & & & & & \\
\to & F_{1,2}'(-1) & {-\phi_{1,2}'\atop\displaystyle\to}  & F_{1,1}'(-1)
& {-\phi_{1,1}'\atop\displaystyle\to}  & F_{1,0}'(-1) &  &  \\
 & & & & & & & \\
    &                &  \hbox to0in{\hss$\scriptstyle -f'_{1,2}$}\searrow &
     & \hbox to0in{\hss$\scriptstyle -f'_{1,1}$} \searrow &     &  \hbox to0in{\hss$\scriptstyle -f'_{1,0}$}\searrow &   \\
 & & & & & & & \\
\to & F_{0,3}'=0 & {\phi_{0,3}'\atop\displaystyle\to} & F_{0,2}'=0 & {\phi_{0,2}'\atop\displaystyle\to}
& F_{0,1}'=0 & {\phi_{0,1}'\atop\displaystyle\to} & F_{0,0}'  \\
 & & & & & & & \end{array}$$
$$\hbox to\hsize{\hfil Figure 1\hfil}$$

\medskip

\begin{lem}\label{cxlem}
The differentials $\phi'_j$, $j>0$, form a complex and
$\phi'_0\phi'_{1}=0;$ i.e., $\phi'_j\phi'_{j+1}=0$ for all
$j\ge0$.
\end{lem}
{\em Proof.} We have
$$\phi'_0\phi'_{1}(s_{k,i,1}\otimes 1)
=\phi'_0(\phi_{i,1}(s_{k,i,1})\otimes 1)=
\phi_{i,0}(\phi_{i,1}(s_{k,i,1}))\otimes x_0^{m-i},$$ which equals
0 since $\phi_{i,j}\phi_{i,j+1}=0$, and we have
\begin{eqnarray}
\phi'_0\phi'_{1}(s_{k,i,0}\otimes1)
& = & \phi'_0(s_{k,i,0}\otimes x_0 - f_{i,0}(s_{k,i,0})\otimes1)\nonumber\\
& = & \phi_{i,0}(s_{k,i,0})\otimes x_0^{m-i+1} -
\phi_{i-1,0}(f_{i,0}(s_{k,i,0}))\otimes x_0^{m-i+1}\nonumber
\end{eqnarray}
which equals 0 since $\phi_{i,0} = \phi_{i-1,0}f_{i,0}$.
Similarly, for $j>0$, we have:
$$\phi'_j\phi'_{j+1}(s_{k,i,j+1}\otimes1)
=\phi'_j(\phi_{i,j+1}(s_{k,i,j+1})\otimes 1)=
\phi_{i,j}(\phi_{i,j+1}(s_{k,i,j+1}))\otimes 1=0$$ and
$$\phi'_j\phi'_{j+1}(s_{k,i,j}\otimes1)  =
\phi'_j(s_{k,i,j}\otimes x_0 -
f_{i,j}(s_{k,i,j})\otimes1-\phi_{i,j}(s_{k,i,j})\otimes1)$$
$$= \phi_{i,j}(s_{k,i,j})\otimes x_0 - \phi_{i-1,j}(f_{i,j}(s_{k,i,j}))\otimes 1$$
$$- (\phi_{i,j}(s_{k,i,j})\otimes x_0 - f_{i,j-1}(\phi_{i,j}(s_{k,i,j}))\otimes1
-\phi_{i,j-1}\phi_{i,j}(s_{k,i,j})\otimes1)$$ which equals 0 since
$\phi_{i-1,j}f_{i,j} = f_{i,j-1}\phi_{i,j}$.
\prfend

The next result is that this complex gives a resolution. This
means first that generators for $I(Z')$ are given by taking for
each $i$ the generators of $I(Z_i)$ times $x_0^{m-i}$, and second
that syzygies for $I(Z')$ are of two types. One type comprises the
obvious syzygies coming from each of the ideals $I(Z_i)$
individually, which give syzygies of $I(Z')$ based on the fact
that $x_0^{m-i}I(Z_i)_t\subset I(Z')_{t+m-i}$. The other type
comprises syzygies between elements of $I(Z_{i+1})$ and elements
of $I(Z_i)$. (These are easy to see too. Not only do we have
$x_0^{m-i}I(Z_i)_t\subset I(Z')_{t+m-i}$, but given a generator
$g$ of $I(Z_{i+1})$ of degree $t$, we get an element
$a=x_0^{m-i-1}g\in I(Z')_{t+m-i-1}$. Now $g\in I(Z_{i+1})_t\subset
I(Z_i)_t$, so we also get an element $b=x_0^{m-i}g\in
I(Z_i)_{t+m-i}\subset I(Z')_{t+m-i}$, and of course $x_0a=b$,
which gives a syzygy of the second type.)

\begin{lem}\label{reslem}
The complex $F'_\bullet$ gives a resolution of $I(Z')$.
\end{lem}
{\em Proof.} We must check that the image of $\phi'_0$ is $I(Z')$
and, for all $j\ge1,$ that the image of $\phi'_j$ is the kernel of
$\phi'_{j-1}$. Let $f\in I(Z')$. We may write $f=x_0g+h$, where no
term of $h$ is divisible by $x_0$. By restricting to the
hyperplane $x_0=0$, we see that $h\in I(Z_m)\subset I(Z')$, hence
that $x_0g\in I(Z')$, and so $g\in I(Z'_{m-1})$. If $m=1$, then
$g\in I(Z'_{0})=R'$, so $f\in x_0^1I(Z_{0})+I(Z_1)$. It follows
that $I(Z')=x_0I(Z_0)+I(Z_1)$. If $m>1$, it follows by induction
on $m$ that $I(Z')=x_0^mI(Z_0)+\cdots+x_0^1I(Z_{m-1})+I(Z_m)$, and
since the image of $\phi_{i,0}$ is $I(Z_i)$, it follows from the
definition that $\phi'_0$ maps onto $I(Z')$, as required.

Now suppose $\phi'_0(f)=0$. Say $f\in F'_{0,0}=R'$. Then $\phi'_0(f)=x_0^mf$
so $f=0$. Now induct on $i$; say $f\in F'_{0,0}\oplus\cdots\oplus
F'_{i,0}$. We can write $f=g+a$, where $g\in
F'_{0,0}\oplus\cdots\oplus F'_{i-1,0}$ and
$a\in F'_{i,0}$. Write $a=c+x_0d$, where $c$ is the sum
of the terms of $a$ not involving $x_0$. Then
$0=\phi'_0(f)=\phi'_0(g)+x_0^{m-i+1}\phi'_{i,0}(d)+x_0^{m-i}\phi'_{i,0}(c)$.
If $\phi'_{i,0}(c)\ne0$, then
$x_0^{m-i+1}$ divides $\phi'_0(g)+x_0^{m-i+1}\phi'_{i,0}(d)$
but not any term of $x_0^{m-i}\phi'_{i,0}(c)$ since $c$ and hence $\phi'_{i,0}(c)$
does not involve $x_0$, so we see that
$\phi'_0(g)+x_0^{m-i+1}\phi'_{i,0}(d)=0$ and
$x_0^{m-i}\phi'_{i,0}(c)=0$. Thus in fact
$\phi'_{i,0}(c)=0$, and hence $c=\phi'_{i,1}(c')$ for some $c'\in
F'_{i,1}$. Since we can regard $d$ as being in $F'_{i,0}(-1)$ and
$\phi'_1(d)=x_0d-f'_{i,0}(d)$, it is
enough now to show $g+f'_{i,0}(d)$ is in the image of $\phi'_1$.
This follows by induction, since $g+f'_{i,0}(d)\in
F'_{0,0}\oplus\cdots\oplus F'_{i-1,0}$ and, by Lemma \ref{cxlem},
$\phi'_0\phi'_1(d)=0$, so $\phi'_0f'_{i,0}(d)=\phi'_0(\mu_i(d))$,
hence $\phi'_0(g+f'_{i,0}(d))=\phi'_0(g+\mu_i(d))=\phi'_0(f)=0$.

Next, suppose $\phi'_j(f)=0$ for some $f$ and some $j\geq 1$.
Inducting on $i$, say $f\in(F'_{0,j})\oplus(F'_{1,j}\oplus F'_{1,j-1}(-1))\oplus\cdots\oplus
(F'_{i,j}\oplus F'_{i,j-1}(-1))$. The case $i=0$ is immediate, so say $i\ge1$.
We can write $f=g+a+b$,
where $g\in (F'_{0,j})\oplus(F'_{1,j}\oplus F'_{1,j-1}(-1))\oplus\cdots\oplus (F'_{i-1,j}\oplus F'_{i-1,j-1}(-1))$,
$a\in F'_{i,j}$ and $b\in F'_{i,j-1}(-1)$. Now
$0=\phi'_j(f)=\phi'_j(g)-f'_{i,j-1}(b)-\phi'_{i,j-1}(b)+\phi'_{i,j}(a)+x_0b$,
when $j>1$, and
$0=\phi'_1(f)=\phi'_1(g)-f'_{i,0}(b)+\phi'_{i,1}(a)+x_0b$ for $j=1$. But
$\phi'_j(g)-f'_{i,j-1}(b)-\phi'_{i,j-1}(b)\in 
(F'_{0,j})\oplus(F'_{1,j}\oplus F'_{1,j-1}(-1))\oplus\cdots\oplus (F'_{i-1,j}\oplus
F'_{i-1,j-1}(-1))$, for $j>1$, $\phi'_1(g)-f'_{i,0}(b)\in
F'_{0,0}\oplus\cdots\oplus F'_{i-1,0}$ for $j=1$ and
$\phi'_{i,j}(a)+x_0b\in F'_{i,j-1}$, so each is 0.

Denote by $c$ the sum of all terms of $a$ not divisible by $x_0$. Then
$a=c+x_0d$ for some $d\in F'_{i,j}(-1)$.
As above we must have $\phi'_{i,j}(c)=0$ (hence
$c=\phi'_{i,j+1}(e)$ for some $e\in F'_{i,j+1}$) and
$\phi'_{j+1}(d)=x_0d-f'_{i,j}(d)-\phi'_{i,j}(d)
=x_0d-f'_{i,j}(d)+b$ (since
$x_0b=-\phi'_{i,j}(a)=-\phi'_{i,j}(c+x_0d)=-\phi'_{i,j}(x_0d)$), so
$\phi'_{j+1}(e+d)=a+b-f'_{i,j}(d)$. Thus it is enough to
show that $f-\phi'_{j+1}(e+d)=g+f'_{i,j}(d)$ is in the image of
$\phi'_{j+1}$. But $g+f'_{i,j}(d)\in 
(F'_{0,j})\oplus(F'_{1,j}\oplus F'_{1,j-1}(-1))\oplus\cdots\oplus (F'_{i-1,j}\oplus
F'_{i-1,j-1}(-1))$ and $\phi'_j(g+f'_{i,j}(d))=\phi'_j(g)-f'_{i,j-1}(b)=0$, 
so this follows by induction.
\prfend

\begin{cor}\label{mincor}
If $I(Z_{i+1})\subset R_1I(Z_{i})$ for each $i$,
then the maps $f_{i,j}$ can be chosen so that
the resolution $F'_\bullet$ is minimal.
\end{cor}
{\em Proof.} By Corollary \ref{fundcor}, we may assume
$f_{i+1,j}:F_{i+1,j}\to R_1F_{i,j}$, for all $i$ and $j$. Since
the resolutions $F_{i,\bullet}$ are minimal, we know that the
matrix for each map $\phi_{i,j}$ has entries in $R_1$, and we now
know the same is true for each $f_{i+1,j}$. The same now follows
for each $\phi'_j$ by an inspection of the definition of
$\phi'_j$. \prfend

We can now prove Theorem \ref{Pointhm}.
\vskip\baselineskip

\noindent{\em Proof of  Theorem \ref{Pointhm}}. Since
$F'_0=\bigoplus_{i=0}^m F'_{i,0}$ and, for $j\ge1$,
$F'_j=\bigoplus_{i=1}^m (F_{i,j}'\oplus F_{i,j-1}'(-1))$, then,
accounting for shifts, we see that the Poincar\'e polynomial
$P(Z')$ is given by the sum
$T^m+XT^{m}P(Z_1)+T^{m-1}P(Z_1)+\cdots+ XTP(Z_m)+P(Z_m)$, which
simplifies to $T^m+(1+XT)(T^{m-1}P(Z_1)+\cdots+P(Z_m))$, as
claimed. \prfend

\section{Applications}\label{appls}

To apply our results to obtain minimal resolutions, we need to
verify the condition $I(Z_{i+1})\subset R_1I(Z_{i})$ of Corollary
\ref{mincor}. We first do this when $\hbox{char}(K)$ is either 0 or
is sufficiently large, then in various additional situations,
such as the case of monomial ideals.

\begin{prop}\label{char0prop}
Given $p_i\in \pr d$ and fat points
$Z=m_1'p_1+\ldots+m_r'p_r$ and $Y=m_1p_1+\ldots+m_rp_r$ with
$m_i'>m_i$ for all $i$ whenever $m_i'>0$, assume
$\hbox{char}(K)$ is either 0 or bigger than the degree
of each generator in a minimal set of homogeneous generators
of $I(Z)$. Then $I(Z)\subset R_1I(Y)$.
\end{prop}
{\em Proof.} This follows easily using Euler's identity,
that $\delta F=\sum_ix_i\partial F/\partial x_i$
for any homogeneous form $F$ with either $\hbox{char}(K)=0$
or $\hbox{char}(K)>\delta$,
where $\delta$ is the degree of $F$. \prfend

\begin{prop}\label{monomprop}
Let $Y$ and $Z$ be as in Proposition \ref{char0prop}. Assume that the
points $p_i$ are located at coordinate vertices of
$\pr {d-1}$. Then $I(Z)\subset R_1I(Y)$.
\end{prop}
{\em Proof.} The ideals $I(Y)$ and $I(Z)$ are generated by
monomials in this case, and we may assume that the variables
are indexed so that $x_j$ vanishes at $p_i$ for all $i\ne j$.
Now assume there is a monomial $f=x_1^{n_1}\cdots
x_d^{n_d}\in (I(Z)\setminus R_1I(Y))\subset I(Y)$. Since $f\in
I(Z)$, we know $m_i'\le (n_1+\cdots+n_d)-n_i$ for every $i$, but
$f\in I(Y)\setminus R_1I(Y)$, so $m_i=(n_1+\cdots+n_d)-n_i$ for
some $i$. However, $m_i'>m_i$, so this is impossible. \prfend

We also have the following bootstrapping result:

\begin{prop}\label{bootprop}
Let $p_1,\ldots,p_r\in L$, where $L\subset \pr {d}$ is a
proper linear subspace of $\pr {d}$. Let $R$ be the
homogeneous coordinate ring for $L$, and $R'$ that for
$\pr {d}$. Given positive integers $m_1,\ldots,m_r$, let
$Z=m_1p_1+\cdots+m_rp_r\subset L$ be the  fat point subscheme of
$L$, and let $Z'=m_1p_1+\cdots+m_rp_r\subset \pr {d}$ be the
fat point subscheme of $\pr {d}$ specified by the same
multiplicities. If $I(Z_{i+1})\subset R_1I(Z_{i})$ holds for all
$i$, then so does $I(Z'_{i+1})\subset R_1I(Z'_{i})$.
\end{prop}
{\em Proof.} It is enough by induction to prove this in the case that
$L$ is a hyperplane. But by Lemma \ref{reslem}, we have
that $I(Z'_i)=\sum_{0\le j\le i}x_0^{i-j}I(Z_j)R'$, for all $i$.
Thus $I(Z'_{i+1})=\sum_{0\le j\le i+1}x_0^{i+1-j}I(Z_j)R'
=x_0(\sum_{0\le j\le i}x_0^{i-j}I(Z_j)R')+I(Z_{i+1})R'$, and this is a subset
of $x_0(\sum_{0\le j\le i}x_0^{i-j}I(Z_j)R')+R_1(I(Z_i)R')
\subset R'_1I(Z'_i)$.
\prfend

\begin{example}\label{3pts}\rm
Resolutions of ideals for fat point subschemes supported at up to $d+1$
general points of $\pr d$ are known in various cases (\cite {fatabbi}, \cite {fatlor}, \cite
{franc}, \cite {valla}). Proposition \ref{monomprop}, and Corollary \ref{mincor}
(or Theorem \ref{Pointhm} for just the Betti numbers),
reduce the problem of determining resolutions of fat point subschemes with support
at up to $d+1$ general points of $\pr d$ to cases in which the
support spans the entire projective space. For example,
to determine the resolution for $2p_1+2p_2+p_3$ for general points $p_i\in \pr d$
with $d>2$, it is enough to determine the resolutions
of $2p_1+2p_2+p_3\subset \pr 2$ and $p_1+p_2\subset \pr 2$,
and to do $p_1+p_2\subset \pr 2$ it suffices to do $p_1+p_2\subset \pr 1$.
Since, in fact, resolutions for ideals of fat points with support at
3 general points of $\pr 2$ are known (\cite{marvi}), our results
as a consequence give the resolution and Betti numbers for ideals
of fat points supported at any 3 general points in projective space of any dimension.
This generalizes the result for two points (\cite {fatlor}, \cite {valla}).
\end{example}

\begin{example}\label{collinearpts} \rm
Another way to generalize the known resolution of fat points with
support at two points is to consider supports consisting of
collinear points. Let $p_1,\ldots,p_r\in L$, where $L\subset
\pr {d}$ is a line. Let $Z=m_1p_1+\cdots+m_rp_r\subset L$.
Then $I(Z)=fR$, where $f$ is a polynomial vanishing at each point
$p_i$ to order $m_i$ and $R$ is the homogeneous
coordinate ring of $L$. Since the ideals are principal and the
degree of the generator is the sum of the multiplicities, it is
easy to see that $I(Z_{i+1})\subset R_1I(Z_{i})$ holds for all
$i$. Proposition \ref{bootprop} and Corollary \ref{mincor} now
give the minimal resolution of $I(Z')\subset \pr {d}$,
and Theorem \ref{Pointhm} gives the Betti numbers,
extending the result of \cite{freeres} for $d=2$.
\end{example}

\begin{example}\label{ptsinp2} \rm
Various facts are known for resolutions of points in $\pr 2$; our
results thus extend these to higher dimension, at least in
characteristic 0. For example, \cite{marvi} works out the
resolution of the ideal of fat points with support on an irreducible
conic in $\pr 2$, while \cite{freeres} determines the Betti
numbers for the case of any plane conic, irreducible or not, and \cite{FHH}
determines the Betti numbers for any fat point subscheme with
support at up to 8 general points of $\pr2$. Thus our results give
the resolution for the ideal of fat points whose support lies on an
irreducible conic curve in any projective space
(since an irreducible degree 2 curve is contained in a plane), 
and they give the Betti numbers when the support either consists of up to 8
general points in a plane or lies in any conic in a plane, for a
plane in any projective space. There are also many additional
examples of sets of points $p_1,\ldots,p_r$ contained in configurations
of lines in the plane for which the graded Betti numbers for
both $p_1+\cdots+p_r$ and $2p_1+\cdots+2p_r$  are known (see \cite{GMS}).
Our results thus give the graded Betti numbers for these examples
regarded as subschemes of $\pr n$, by regarding $\pr2$ as a linear
subspace of $\pr n$.
\end{example}

\begin{remark}\label{finalrem} \rm
We close with a remark about an additional situation
in which our criterion for minimality will hold.
Consider a fat point subscheme $Z\subset \pr d$.
For each $i$, let $D_i$ ($d_i$, resp.) be the degree of
the generator of maximal (resp., minimal) degree in a
minimal set of homogeneous generators for $I(Z_i)$.
Since $I(Z_{i+1})\subset I(Z_i)$ and
$R_1I(Z_i)_t=I(Z_i)_{t+1}$ for $t\ge D_i$,
it is clear that the condition $I(Z_{i+1})\subset
R_1I(Z_i)$ holds if the degrees of the generators of
$I(Z_{i+1})$ are shifted enough with respect to those of $I(Z_i)$
(in particular, if $d_{i+1}>D_i$ for each $i>0$).
This occurs, for example, if, for each $i$, the fat
points in $Z_{i+1}$ of multiplicity 1 are general and if
there are enough of them. More explicitly,
let $Z_1$ consist of $r_0$ simple points.
Let $Z_2$ include the same
points as does $Z_1$, but take these points with multiplicity 2,
and add on $r_1$ additional general simple points.
Continue in this way, defining $Z_i$ and $r_i$.
The condition $I(Z_{i+1})\subset R_1I(Z_i)$ holds for all $i$, if,
for example, $r_i\ge {D_i+d\choose d}$ for all $i$, since
$I(Z_{i+1})$ has no elements of degree less than $D_i+1$.
\end{remark}

\end{document}